\input amstex
\magnification=\magstep1 
\baselineskip=13pt
\documentstyle{amsppt}
\vsize=8.7truein \CenteredTagsOnSplits \NoRunningHeads
 \def\per{\operatorname{per}}
 \def\vl{\operatorname{vol}}
 \def\PP{\bold{Pr}\thinspace}
 \def\EE{\bold{E}\thinspace}
 \def\HH{\bold{H}}
\def\xx{\bold{x}}
\def\yy{\bold{y}}
\def\sss{\bold{s}}
\def\ttt{\bold{t}}

 \topmatter

\title  Matrices with prescribed row and column sums \endtitle
\author Alexander Barvinok  \endauthor
\address Department of Mathematics, University of Michigan, Ann Arbor,
MI 48109-1043, USA \endaddress
\email barvinok$\@$umich.edu \endemail
\date October 2010 \enddate
\thanks  This research  was partially supported by NSF Grant DMS 0856640 and
by a United States - Israel BSF grant 2006377.
\endthanks 
\abstract 
 This is a survey of the recent progress and open questions on the structure of the 
sets of 0-1 and non-negative integer matrices with prescribed row and column sums. 
We discuss cardinality estimates, the structure of a random matrix from the set, discrete versions of the Brunn-Minkowski inequality and the statistical dependence between row and column sums.
\endabstract
\endtopmatter

\document

\head 1. Introduction  \endhead

Let $R=\left(r_1, \ldots, r_m \right)$ and $C=\left(c_1, \ldots, c_n \right)$
be positive integer vectors such that 
$$r_1 + \ldots + r_m = c_1 + \ldots + c_n =N. \tag1.1$$
We consider the set $A_0(R, C)$ of all $m \times n$ matrices $D=\left(d_{ij}\right)$ with 0-1 entries, row 
sums $R$ and column sums $C$:
$$\split A_0(R, C)=\Biggl\{ D=\left(d_{ij}\right): \qquad &\sum_{j=1}^n d_{ij} =r_i \quad \text{for} \quad i=1, \ldots, m\\
&\sum_{i=1}^m d_{ij} =c_j \quad \text{for} \quad j=1, \ldots, n \\
&d_{ij} \in \{0, 1\} \Biggr\}. \endsplit$$

We also consider the set $A_+(R, C)$ of non-negative integer $m \times n$ matrices with row sums $R$ and column 
sums $C$:
$$\split A_+(R, C)=\Biggl\{ D=\left(d_{ij}\right): \qquad &\sum_{j=1}^n d_{ij} =r_i \quad \text{for} \quad i=1, \ldots, m\\
&\sum_{i=1}^m d_{ij} =c_j \quad \text{for} \quad j=1, \ldots, n \\
&d_{ij} \in {\Bbb Z}_+ \Biggr\}. \endsplit$$

Vectors $R$ and $C$ are called {\it margins} of matrices from $A_0(R, C)$ and $A_+(R, C)$. We reserve notation $N$ for the sums 
of the coordinates of $R$ and $C$ in (1.1) and write $|R|=|C|=N$.

While the set $A_+(R, C)$ is non-empty as long as the balance condition (1.1) is satisfied,  a result of Gale and 
Ryser (see, for example, Section 6.2 of \cite{BR91}) provides a necessary and sufficient criterion for set 
$A_0(R, C)$ to be non-empty. Let us assume that 
$$\split &m \ \geq \ c_1 \ \geq \ c_2 \ \geq \ \ldots \ \geq \ c_n \ \geq \ 0 \quad \text{and that} \\
& n \ \geq \ r_i \ \geq \ 0
\quad \text{for} \quad i=1, \ldots, n. \endsplit$$
 Set $A_0(R,C)$ is not empty 
if and only if (1.1) holds and 
$$\sum_{i=1}^m \min\left\{r_i, k \right\} \ \geq \ \sum_{j=1}^k c_j \quad \text{for} \quad k=1, \ldots, n.$$
Assuming that $A_0(R, C) \ne \emptyset$, we are interested in the following questions:
\bigskip
$\bullet$ What is the cardinality $\left|A_0(R, C)\right|$ of $A_0(R, C)$ and the cardinality $\left|A_+(R, C)\right|$ of
$A_+(R, C)$?
\medskip
$\bullet$ Let us us consider $A_0(R, C)$ and $A_+(R, C)$ as finite probability spaces with the uniform measure. What a random matrix 
$D \in A_0(R, C)$ and a random matrix $D \in A_+(R, C)$ are likely to look like?
\bigskip
The paper is organized as follows.

In Section 2 we estimate of $\left|A_0(R, C)\right|$ within an $(mn)^{O(m+n)}$ factor and in Section 3 we estimate
$\left| A_+(R, C)\right|$ within an $N^{O(m+n)}$ factor. In all but very sparse cases this way
we obtain asymptotically exact estimates of $\ln \left| A_0(R, C)\right|$ and $\ln \left| A_+(R, C)\right|$ respectively.
The estimate of Section 2 is based on 
a representation of $\left| A_0(R, C)\right|$ as the permanent of a certain $mn \times mn$ matrix of 0's and 1's, while 
the estimate of Section 3 is based on a representation of $\left| A_+(R, C) \right|$ as the expectation of the 
permanent of a certain $N \times N$ random matrix with exponentially distributed entries. 
In the proofs, the crucial role is played by 
the van der Waerden inequality for permanents of doubly stochastic matrices. The cardinality estimates are obtained as solutions 
to simple convex optimization problems and hence are efficiently computable, although they 
cannot be expressed by a ``closed formula'' in the margins $(R, C)$. 
Our method is sufficiently robust as the same approach can be 
applied to estimate the cardinality of the set of matrices with prescribed margins {\it and} with 0's in prescribed positions.

In Sections 4 and 5 we discuss some consequences of the formulas obtained in Sections 2 and 3. In particular,
in Section 4, we show that the numbers $\left| A_0(R, C)\right|$ and $\left| A_+(R, C)\right|$ are both 
approximately log-concave as functions of the margins $(R, C)$. We note an open question whether these numbers 
are genuinely log-concave and give some, admittedly weak, evidence that it may be the case. In Section 5, we 
discuss statistical dependence between row and column sums. Namely, we consider finite probability spaces of $m \times n$
non-negative integer or 0-1 matrices with the total sum $N$ of entries and two events in those spaces: event 
${\Cal R}$ consisting of the matrices with row sums $R$ and event ${\Cal C}$ consisting of the matrices with 
column sums $C$. It turns out that 0-1 and non-negative integer matrices exhibit opposite types of behavior.
Assuming that the margins $R$ and $C$ are sufficiently far away from sparse and uniform, we show that 
for 0-1 matrices the events ${\Cal R}$ and ${\Cal C}$ repel each other (events ${\Cal R}$ and ${\Cal C}$ 
are negatively correlated) while for 
non-negative integer matrices they attract each other (the events are positively correlated). 

In Section 6, we discuss what random matrices $D \in A_0(R, C)$ and $D \in A_+(R, C)$ look like.
We show that in many respects, a random matrix $D \in A_0(R, C)$ behaves like an $m \times n$ matrix $X$ of
independent Bernoulli random variables such that $\EE X=Z_0$ where $Z_0$ is a certain matrix, called 
the maximum entropy matrix, with row sums $R$,
column sums $C$ and entries between 0 and 1. It turns out that $Z_0$ is the solution to an optimization problem,
which is convex dual to the optimization problem of Section 2 used to estimate $\left| A_0(R, C)\right|$.
On the other hand, a random matrix $D \in A_+(R, C)$ in many respects behaves like an $m \times n$ matrix 
$X$ of independent geometric random variables such that $\EE X=Z_+$ where $Z_+$ is a certain matrix, also called the 
maximum entropy matrix, with row sums $R$, column sums $C$ and non-negative entries. It turns out that $Z_+$ is 
the solution to an optimization problem which is convex dual to the optimization problem of Section 3 used to estimate
$\left| A_+(R, C)\right|$. It follows that in various natural metrics matrices $D \in A_0(R, C)$ concentrate about 
$Z_0$ while matrices $D \in A_+(R, C)$ concentrate about $Z_+$. We note some open questions on whether 
individual entries of  random $D \in A_0(R, C)$ and random $D \in A_+(R, C)$ are asymptotically Bernoulli, respectively geometric,
with the expectations read off from $Z_0$ and $Z_+$.

In Section 7, we discuss asymptotically exact formulas for $\left| A_0(R, C) \right|$ and \break $|A_+(R, C)|$. Those formulas 
are established under essentially more restrictive conditions than cruder estimates of Sections 2 and 3. We assume that the entries of the maximum entropy matrices $Z_0$ and $Z_+$ are within a constant factor,
fixed in advance, of each other. Recall that matrices $Z_0$ and $Z_+$ characterize the 
typical behavior of random matrices $D \in A_0(R, C)$ and $D \in A_+(R, C)$ respectively.
In the case of 0-1 matrices our condition basically means that the margins $(R, C)$ lie 
sufficiently deep inside the region defined by the Gale-Ryser inequalities. As the margins approach the boundary, the 
number $\left| A_0(R, C)\right|$ gets volatile and hence cannot be expressed by an analytic formula like the 
one described in Section 7. The situation with non-negative integer matrices is less clear. It is plausible that the 
number $\left| A_+(R, C)\right|$ experiences some volatility when some entries of $Z_+$ become abnormally large, but 
we don't have a proof of that happenning.

In Section 8, we mention some possible ramifications, such as enumeration of higher-order tensors and graphs 
with given degree sequences.

The paper is a survey and although we don't provide complete proofs, we often sketch main ideas of our approach.

\head 2. The logarithmic asymptotic for the number of 0-1 matrices  \endhead

The following result is proven in \cite{Ba10a}.

\proclaim{(2.1) Theorem} Given positive integer vectors 
$$R=\left(r_1, \ldots, r_m\right) \quad \text{and} \quad  C=\left(c_1, \ldots, c_n \right),$$
 let us define the function
$$\split &F_0(\xx, \yy)=\left(\prod_{i=1}^m x_i^{-r_i} \right) \left(\prod_{j=1}^n y_j^{-c_j}\right) \left(\prod_{i,j} \left(1+x_i y_j \right)
\right)\\
&\quad \text{for} \quad \xx=\left(x_1, \ldots, x_m \right) \quad \text{and} \quad \yy=\left(y_1, \ldots, y_n \right)
\endsplit$$
and let
$$\alpha_0(R, C)=\inf \Sb x_1, \ldots, x_m >0 \\ y_1, \ldots, y_n > 0 \endSb F_0(\xx, \yy).$$
Then the number $A_0(R, C)$ of $m \times n$ zero-one matrices with row sums $R$ and column sums $C$ satisfies
$$\alpha_0(R, C) \ \geq \ |A_0(R, C)| \ \geq \ {(mn)! \over (mn)^{mn}} 
\left( \prod_{i=1}^m {\left(n-r_i\right)^{n-r_i} \over \left(n-r_i\right)!} \right) \left( \prod_{j=1}^n {c_j^{c_j} \over c_j!}\right)
\alpha_0(R, C).$$
\endproclaim
Using Stirling's formula, 
$${s! \over s^s}=\sqrt{2 \pi s} e^{-s} \left(1 +O\left({1 \over s}\right)\right),$$
one can notice that the ratio between the upper and lower bounds is $(mn)^{O(m+n)}$. Indeed, the 
``$e^{-s}$'' terms cancel each other out, since 
$$e^{-mn} \left(\prod_{i=1}^m e^{n-r_i}\right)\left( \prod_{j=1}^n e^{c_j}\right)=1.$$
Thus, for sufficiently dense 0-1 matrices, where we have $|A_0(R, C)| =2^{\Omega(mn)}$, 
we have an asymptotically exact formula
$$\ln |A_0(R, C)|\ \approx \ \ln \alpha_0(R, C) \quad \text{as} \quad m, n \longrightarrow +\infty.$$

\subhead (2.2) A convex version of the optimization problem \endsubhead 
Let us substitute 
$$x_i =e^{s_i} \quad \text{for} \quad i=1, \ldots, m \quad \text{and} \quad 
y_j=e^{t_j} \quad \text{for} \quad j=1, \ldots, n$$
in $F_0(\xx, \yy)$.
Denoting
$$\split &G_0(\sss, \ttt)=-\sum_{i=1}^m r_i s_i -\sum_{j=1}^n t_j c_j + \sum_{i,j} \ln \left(1+e^{s_i +t_j}\right)\\
&\quad \text{for} \quad \sss=\left(s_1, \ldots, s_m\right) \quad \text{and} \quad \ttt=\left(t_1, \ldots, t_n \right),
\endsplit \tag2.2.1$$
we obtain
$$\ln \alpha_0(R, C) =\inf \Sb s_1, \ldots, s_m \\ t_1, \ldots, t_n \endSb G_0(\sss, \ttt).$$
We observe that $G_0(\sss, \ttt)$ is a convex function on ${\Bbb R}^{m+n}$. In particular, one can compute 
the infimum of $G_0$ efficiently by using interior point methods, see, for example, \cite{NN94}.

\subhead (2.3) Sketch of proof of Theorem 2.1 \endsubhead
The upper bound for $\left|A_0(R,C)\right|$ is immediate: it follows from the 
expansion
$$\prod_{ij} \left(1+ x_i y_j\right) = \sum_{R, C} |A_0(R, C)| \xx^R \yy^C,$$
where 
$$\xx^R=x_1^{r_1} \cdots x_m^{r_m} \quad \text{and} \quad \yy^C=y_1^{c_1} \cdots y_n^{c_n}$$ 
for $R=\left(r_1, \ldots, r_m \right)$ and $C=\left(c_1, \ldots, c_n \right)$ and the sum is taken over 
all pairs of non-negative integer vectors $R=\left(r_1, \ldots, r_m\right)$ and $C=\left(c_1, \ldots, c_n \right)$ 
such that $r_1 + \ldots +r_m =c_1 + \ldots +c_n \leq mn$.

To prove the lower bound, we express $|A_0(R, C)|$ as the permanent of an $mn \times mn$ matrix. 
Recall that the {\it permanent} of a $k \times k$ matrix $B=\left(b_{ij}\right)$ 
is defined by 
$$\per B =\sum_{\sigma \in S_k} \prod_{i=1}^k b_{i \sigma(i)},$$
where the sum is taken over the symmetric group $S_k$ of all permutations $\sigma$ of the set 
$\{1, \ldots, k\}$, see, for example, Chapter 11 of \cite{LW01}. 
One can show, see \cite{Bar10a} for details, that 
$$|A_0(R, C)| = \left( \prod_{i=1}^m {1 \over (n-r_i)!} \right) \left(\prod_{j=1}^n {1 \over c_j!} \right)
\per B, \tag2.3.1$$
where $B$ is the $mn \times mn$ matrix of the following structure:
\bigskip
the rows of $B$ are split into distinct $m+n$ blocks, the $m$ blocks of type I 
 having $n-r_1, \ldots, n-r_m$ rows respectively and $n$ blocks of type II having $c_1, \ldots, c_n$
 rows respectively; 
 \smallskip
the columns of $B$ are split into $m$ distinct blocks of $n$ columns each;
\smallskip
for $i=1, \ldots, m$, the entry of $B$ that lies in a row from the $i$-th block of rows of type I and 
a column from the $i$-th block of columns is equal to 1;
\smallskip
for $i=1, \ldots, m$ and $j=1, \ldots, n$, the entry of $B$ that lies in a row from the $j$-th block of 
rows of type II and the $j$-th column from the $i$-th block of columns is equal to 1;
\smallskip
all other entries of $B$ are $0$.
\bigskip
Suppose that the infimum of function $G_0(\sss, \ttt)$ defined by (2.2.1) is attained at 
a particular point $\sss=\left(s_1, \ldots, s_m \right)$ and 
$\ttt =\left(t_1, \ldots, t_n \right)$ (the case when the infimum is not attained is handled by an approximation 
argument). Let $x_i =\exp\left\{s_i\right\}$ for $i=1, \ldots, m$ and 
$y_j=\exp\left\{t_j\right\}$ for $j=1, \ldots, n$.

Setting the gradient of $G_0(\sss, \ttt)$ to 0, we obtain
$$\split &\sum_{j=1}^n {x_i y_j \over 1+ x_i y_j} =r_i 
\quad \text{for} \quad i=1, \ldots, m\\
&\sum_{i=1}^m {x_i y_j \over 1+x_i y_j} =c_j \quad \text{for} \quad j=1, \ldots, n.\endsplit \tag2.3.2 $$
Let us consider a matrix $B'$ obtained from matrix $B$ as follows: 
\bigskip
for $i=1, \ldots, m$ we multiply every row of $B$ in the $i$-th block of type I by 
$${1 \over  x_i (n-r_i)};$$

for $j=1, \ldots, n$, we multiply every row of $B$ in the $j$-th block of type II by 
$${y_j \over c_j};$$

for $i=1, \ldots, m$ and $j=1, \ldots, n$ we multiply the $j$-th column in the $i$-th block of 
columns of $B$ by 
$${x_i \over 1+x_i y_j}.$$
\bigskip
Then 
$$\per B = \left(\prod_{i=1}^m x_i^{-r_i} \left(n-r_i\right)^{n-r_i} \right) 
\left(\prod_{j=1}^n y_j^{-c_j} c_j^{c_j} \right) \left(\prod_{ij} \left(1+x_i y_j \right) \right)
\per B'.$$
On the other hand, equations (2.3.2) imply that the row and column sums of $B'$ are equal to 1, 
that is, $B'$ is {\it doubly stochastic}. Applying the van der Waerden bound for permanents 
of doubly stochastic matrices, see, for example, Chapter 12 of \cite{LW01}, we conclude that 
$$\per B' \ \geq \ {(mn)! \over (mn)^{mn}},$$
which, together with (2.3.1) completes the proof.
{\hfill \hfill \hfill} \qed

One can prove a version of Theorem 2.1 for 0-1 matrices with prescribed row and 
column sums {\it and} prescribed zeros in some positions.

\head 3. The logarithmic asymptotics for the number of non-negative integer matrices  \endhead 

The following result is proven in \cite{Ba09}.
\proclaim{(3.1) Theorem} Let $R=\left(r_1, \ldots, r_m\right)$ and $C=\left(c_1, \ldots, c_n \right)$
be positive integer vectors such that $r_1 + \ldots + r_m =c_1 + \ldots +c_n=N$.
Let us define a function 
$$\split &F_+(\xx, \yy)=\left( \prod_{i=1}^m x_i^{-r_i} \right) \left( \prod_{j=1}^n y_j^{-c_j} \right) 
\left( \prod_{ij} {1 \over 1- x_i y_j} \right) \\
&\quad \text{for} \quad \xx=\left(x_1, \ldots, x_m \right) \quad \text{and} \quad 
\yy=\left(y_1, \ldots, y_n \right). \endsplit$$
Then $F_+(\xx, \yy)$ attains its minimum
$$\alpha_+(R, C)=\min \Sb 0 < x_1, \ldots, x_m < 1 \\ 0 < y_1, \ldots, y_n < 1 \endSb F_+(\xx, \yy)$$
on the open cube $0< x_i, y_j < 1$ and for the number $|A_+(R, C)|$ of non-negative 
integer $m \times n$ matrices with row sums $R$ and column sums $C$, we have 
$$\alpha_+(R, C) \ \geq \ | A_+(R, C)| \ \geq \ N^{-\gamma(m+n)} \alpha_+(R, C),$$
where $\gamma > 0$ is an absolute constant.
\endproclaim
For sufficiently dense matrices, where 
$$\min_{i=1, \ldots, m} r_i = \Omega(n) \quad \text{and} \quad \min_{j=1, \ldots, n} c_j =\Omega(m)$$
we have $|A_+(R, C)| =(N/mn)^{\Omega(mn)}$ and hence we obtain an asymptotically exact formula
$$\ln | A_+(R, C)| \ \approx \ \ln \alpha_+(R, C) \quad \text{as} \quad m, n \longrightarrow +\infty.$$

\subhead (3.2) A convex version of the optimization problem \endsubhead 
Let us substitute 
$$x_i =e^{-s_i} \quad \text{for} \quad i=1, \ldots, m \quad \text{and} \quad 
y_j = e^{-t_j} \quad \text{for} \quad j=1, \ldots, n$$
in $F_+(\xx, \yy)$.
Denoting 
$$\split &G_+(\sss, \ttt)=\sum_{i=1}^m r_i s_i + \sum_{j=1}^n t_j c_j - \sum_{i,j} \ln \left(1 - e^{-s_i -t_j}\right) 
\\ 
&\quad \text{for} \quad \sss=\left(s_1, \ldots, s_m \right) \quad \text{and} \quad 
\ttt=\left(t_1, \ldots, t_n \right), \endsplit \tag3.2.1$$
we obtain 
$$\ln \alpha_+(R, C)= \min \Sb s_1, \ldots, s_m > 0 \\ t_1, \ldots t_n > 0 \endSb G_+(\sss, \ttt).$$
We observe that $G_+(\sss, \ttt)$ is a convex function on ${\Bbb R}^{m+n}$. In particular, 
one can compute the minimum of $G_+$ efficiently by using interior point methods \cite{NN94}.

\subhead (3.3) Sketch of proof of Theorem 3.1 \endsubhead 
The upper bound for $|A_+(R, C)|$ follows immediately from the expansion 
$$\prod_{ij} {1 \over 1-x_i y_j} =\sum_{R, C} |A_+(R, C)| \xx^R \yy^C \quad 
\text{for} \quad 0 < x_1, \ldots, x_m, y_1, \ldots, y_n <1$$
where 
$$\xx^R = x_1^{r_1} \cdots x_m^{r_m} \quad \text{and} \quad \yy^C =y_1^{c_1} \cdots y_n^{c_n}$$
for $R=\left(r_1, \ldots, r_m\right)$ and $C=\left(c_1, \ldots, c_n\right)$ and the sum is taken over all pairs 
of non-negative integer vectors $R=\left(r_1, \ldots, r_m\right)$ and $C=\left(c_1, \ldots, c_n \right)$ such that 
$r_1 + \ldots + r_m =c_1 + \ldots +c_n$.

To prove the lower bound, we express $|A_+(R, C)|$ as the integral of the permanent of 
an $N \times N$ matrix with variable entries.  For an $m \times n$ matrix $Z=\left(z_{ij}\right)$
we define the $N \times N$ matrix $B(Z)$ as follows:
\bigskip
 the rows of $B(Z)$ are split into 
$m$ distinct blocks of sizes $r_1, \ldots, r_m$ respectively;
\smallskip
the columns of $B(Z)$
are split into $n$ distinct blocks of sizes $c_1, \ldots, c_n$ respectively;
\smallskip
for $i=1, \ldots, m$ and $j=1, \ldots, n$, the entry of $B(Z)$ that lies in a row from the 
$i$-th block of rows and in a column from the $j$-th block of columns is $z_{ij}$.
\bigskip 
Then there is a combinatorial identity
$$\per B(Z) = \left( \prod_{i=1}^m r_i! \right) \left( \prod_{j=1}^n c_j! \right) \sum \Sb D \in A_+(R, C) \\
D=\left(d_{ij}\right) \endSb \prod_{ij} {z_{ij}^{d_{ij}} \over d_{ij}!},$$
cf. \cite{Be74}, which implies that 
$$|A_+(R, C)|=\left( \prod_{i=1}^m {1 \over r_i!} \right) \left(\prod_{j=1}^n {1 \over c_j!}\right) 
\int_{{\Bbb R}^{mn}_+} \per B(Z) \exp\left\{ - \sum_{ij} z_{ij} \right\} \ dZ.$$
Here the integral is taken over the set ${\Bbb R}^{mn}_+$ of $m \times n$ matrices $Z$ with 
positive entries.  Let
$\Delta_{mn-1} \subset {\Bbb R}^{mn}_+$ be the standard $(mn-1)$-dimensional simplex defined by the equation 
$$\sum_{ij} z_{ij}=1.$$
Since $\per B(Z)$ is a homogeneous polynomial in $Z$ of degree $N$, 
we have 
$$|A_+(R, C)|={(N+mn-1)! \over \sqrt{mn}} \left( \prod_{i=1}^m {1 \over r_i!} \right) \left(\prod_{j=1}^n {1 \over c_j!}\right) 
\int_{\Delta_{mn-1}} \per B(Z) \ dZ, \tag3.3.1$$
where $dZ$ is the Lebesgue measure on $\Delta_{mn-1}$ induced from ${\Bbb R}^{mn}$.

Let $\sss=\left(s_1, \ldots, s_m\right)$ and $\ttt=\left(t_1, \ldots, t_n\right)$ be the minimum point of function 
$G_+(\sss, \ttt)$ defined by (3.2.1). Let $x_i=\exp\left\{s_i\right\}$ for $i=1, \ldots, m$ and $y_j=\exp\left\{t_j \right\}$ 
for $j=1, \ldots, m$. Setting the gradient of $G_+(\sss, \ttt)$ to 0, we obtain
$$\split &\sum_{j=1}^n {x_i y_j \over 1-x_i y_j} =r_i \quad \text{for} \quad i=1, \ldots, m \\
&\sum_{i=1}^m {x_i y_j \over 1-x_i y_j} =c_j \quad \text{for} \quad j=1, \ldots, n. \endsplit \tag3.3.2$$
Let us consider the affine subspace $L \subset {\Bbb R}^{mn}$ of $m \times n$ matrices $Z=\left(z_{ij}\right)$
defined by the system of equations 
$$\split &\sum_{j=1}^n x_i y_j z_{ij} ={r_i \over N+mn}\quad \text{for} \quad i=1, \ldots, m \quad \text{and} \\
&\sum_{i=1}^m x_i y_j z_{ij} ={c_j\over N+mn} \quad \text{for} \quad j=1, \ldots, n. \endsplit \tag3.3.3$$
We note that $\dim L=(m-1)(n-1)$.

Suppose that $Z \in \Delta_{mn-1} \cap L$ and consider the corresponding matrix $B(Z)$. If we multiply every row 
in the $i$-th block of rows by $x_i\sqrt{N+mn}/r_i$ and every column in the $j$-th block of columns by $y_j\sqrt{N+mn}/c_j$, by (3.3.3) we obtain a doubly stochastic matrix $B'(Z)$ for which we have 
$\per B'(Z) \geq N!/N^N$ by the van der Waerden inequality. Summarizing,
$$ \aligned \per B(Z) \ \geq \ &{N!  \over N^N (N+mn)^N} \left( \prod_{i=1}^m r_i^{r_i} \right) \left(\prod_{j=1}^n c_j^{c_j} \right) \\ &\quad \times \left(\prod_{i=1}^m x_i^{-r_i} \right) \left( \prod_{j=1}^n y_j^{-c_j} \right) \qquad
 \text{for all} \quad Z \in \Delta_{mn-1} \cap L. \endaligned \tag3.3.4$$
 It remains to show that the intersection $\Delta_{mn-1} \cap L$ is sufficiently large, so that the
 contribution of a neighborhood of the intersection to the integral (3.3.1) is sufficiently large. 
 It follows by (3.3.2)--(3.3.3) that $\Delta_{mn-1} \cap L$ contains matrix $Z=\left(z_{ij}\right)$ where 
 $$z_{ij}={1 \over (N+mn) \left(1-x_i y_j \right)} \quad \text{for all} \quad i, j.$$ 
In \cite{Ba09}, we prove a geometric lemma which states that if $\Delta_{d-1} \subset {\Bbb R}^d_+$ is the 
standard $(d-1)$-dimensional simplex that is the intersection of the affine hyperplane $H$ defined by the
equation
$x_1+\ldots +x_d=1$ and the positive orthant
$x_1>0, \ldots, x_d >0$ and if $L \subset H$ is an affine subspace of codimension $k$ in $H$ such that
$L$ contains a point $a \in \Delta_{d-1}$, $a=\left(\alpha_1, \ldots, \alpha_d \right)$, then for the 
volume of the intersection $\Delta_{d-1} \cap L$ we have the lower bound
$$\vl_{d-k-1}\left( \Delta_{d-1} \cap L\right) \ \geq \ {\gamma \over d! \omega_k} d^{d-{1\over 2}} \alpha_1 \ldots \alpha_d,$$
where 
$$\omega_k ={\pi^{k/2} \over \Gamma(k/2+1)}$$
is the volume of the $k$-dimensional unit ball and $\gamma >0$ is an absolute constant.
Applying this estimate in our situation, we conclude that 
$$\vl_{mn-k}\left(\Delta_{mn-1} \cap L \right) \ \geq \ {1 \over (mn)^{O(m+n)}} {e^{mn} \over (N+mn)^{mn}} \prod_{ij} {1 \over 1-x_i y_j},$$
where $k=m+n-1$ or $k=m+n$ depending whether or not $L$ lies in the affine hyperplane $\sum_{ij}  z_{ij}=1$.
This allows us to obtain a similar bound for the volume of a small neighborhood of the intersection
$\Delta_{mn-1} \cap L$. 
Because $\per B(Z)$ is a homogeneous polynomial in $Z$ of degree $N$, 
 inequality (3.3.4) holds in the $\epsilon$-neighborhood of the intersection $\Delta_{mn-1} \cap L$ for $\epsilon =N^{-O(m+n)}$
 up to an $N^{O(m+n)}$ factor. Using it together with (3.3.1), we complete the proof of Theorem 3.1.
{\hfill \hfill \hfill} \qed 

One can prove a version of Theorem 3.1 for non-negative integer matrices with prescribed row and column sums 
{\it and} with prescribed zeros in some positions.

\head 4. Discrete Brunn - Minkowski  inequalities \endhead

Theorems 2.1 and 3.1 allow us to establish approximate log-concavity of the numbers $A_0(R, C)$ 
and $A_+(R, C)$. 

For a non-negative integer vector $B=\left(b_1, \ldots, b_p \right)$, we denote
$$|B|=\sum_{i=1}^p b_i.$$
\proclaim{(4.1) Theorem} Let $R_1, \ldots, R_p$ be positive integer $m$-vectors and let $C_1, \ldots, C_p$ 
be positive integer $n$-vectors such that $|R_1|=|C_1|, \ldots, |R_p|=|C_p|$. 

Let 
$\beta_1, \ldots, \beta_p \geq 0$
be real numbers such that $\beta_1 + \ldots + \beta_p =1$ and such that $R=\beta_1 R_1 +\ldots +\beta_p R_p$
is a positive integer $m$-vector and $C=\beta_1 C_1 +\ldots + \beta_p C_p$ is a positive integer $n$-vector.
Let $N=|R|=|C|$.

Then for some absolute constant $\gamma >0$ we have 
\roster
\item $$(mn)^{\gamma(m+n)} \left|A_0(R, C)\right|  \ \geq \ \prod_{k=1}^p \left| A_0\left(R_k, C_k \right)\right|^{\beta_k} $$
and
\item
$$N^{\gamma(m+n)} \left|A_+(R, C)\right|  \ \geq \ \prod_{k=1}^p \left| A_+\left(R_k, C_k \right)\right|^{\beta_k}.$$
\endroster
\endproclaim
\demo{Proof} Let us denote function $F_0$ of Theorem 2.1 for the pair $\left(R_k, C_k\right)$ by 
$F_k$ and for the pair $\left(R, C\right)$ just by $F$. Then
$$F(\xx, \yy) = \prod_{k=1}^p F_k^{\beta_k}(\xx, \yy) \tag4.1.1$$
and  hence 
$$\alpha_0(R, C) \ \geq \ \prod_{k=1}^p \left(\alpha_0\left(R_k, C_k \right)\right)^{\beta_k}.$$
Part (1) now follows by Theorem 2.1. 

Similarly, we obtain (4.1.1) if we denote function $F_+$ of Theorem 3.1 for the pair $\left(R_k, C_k\right)$ by $F_k$ and for the pair 
$\left(R, C\right)$ just by $F$.
Hence
$$\alpha_+(R, C) \ \geq \ \prod_{k=1}^p \left(\alpha_+\left(R_k, C_k\right)\right)^{\beta_k}.$$
Part (2) now follows by Theorem 3.1.
{\hfill \hfill \hfill} \qed
\enddemo
Teorem 2.1 implies a more precise estimate 
$${(mn)^{mn} \over (mn)!} \left( \prod_{i=1}^m {\left(n-r_i\right)! \over \left(n-r_i\right)^{n-r_i}} \right)
\left(\prod_{j=1}^n {c_j! \over c_j^{c_j}}\right) 
\left| A_0(R, C)\right| \ \geq \ \prod_{k=1}^p \left|A_0\left(R_k, C_k \right)\right|^{\beta_k},$$
where $R=\left(r_1, \ldots, r_m\right)$ and $C=\left(c_1, \ldots, c_n \right)$.

In \cite{Ba07} a more precise estimate 
$${N^N \over N!} \min\left\{ \prod_{i=1}^m {r_i! \over r_i^{r_i}}, \quad 
\prod_{j=1}^n {c_j! \over c_j^{c_j}} \right\} \left| A_+(R, C)\right| \ \geq \ 
\prod_{k=1}^p \left| A_+\left(R_k, C_k \right)\right|^{\beta_k}$$ 
is proven
under the additional assumption that $|R_k|=|C_k|=N$ for $k=1, \ldots, p$. 

Theorem 4.1 raises a natural question whether stronger inequalities hold. 

\subhead (4.2) Brunn-Minkowski inequalities \endsubhead 

\subsubhead (4.2.1) Question \endsubsubhead Is it true that under the conditions of Theorem 4.1 we have 
$$\left| A_0(R, C)\right| \ \geq \ \prod_{k=1}^p \left| A_0\left(R_k, C_k\right)\right|^{\beta_k} ?$$

\subsubhead (4.2.2) Question \endsubsubhead Is it true that under the conditions of Theorem 4.1 we have 
$$\left| A_+(R, C)\right| \ \geq \ \prod_{k=1}^p \left| A_+\left(R_k, C_k\right)\right|^{\beta_k} ?$$

Should they hold, inequalities of (4.2.1) and (4.2.2) would be natural examples of discrete Brunn-Minkowski
inequalities, see \cite{Ga02} for a survey.

Some known simpler inequalities are consistent with the inequalities of (4.2.1)--(4.2.2). 
Let $X=\left(x_1, \ldots, x_p\right)$ and $Y=\left(y_1, \ldots, y_p \right)$ be non-negative 
integer vectors such that 
$$x_1 \geq x_2 \geq \ldots \geq x_p \quad \text{and} \quad y_1 \geq y_2 \geq \ldots \geq y_p.$$
We say that $X$ {\it dominates} $Y$ if 
$$\sum_{i=1}^k x_i \ \geq \ \sum_{i=1}^k y_i \quad \text{for} \quad k=1, \ldots, p-1 
\quad \text{and} \quad \sum_{i=1}^p x_i = \sum_{i=1}^p y_i.$$ 
Equivalently, $X$ dominates $Y$ if $Y$ is a convex combination of vectors obtained from 
$X$ by permutations of coordinates. 

One can show that 
$$\left| A_0(R, C)\right| \ \geq \ \left| A_0\left(R', C'\right)\right| \quad 
\text{and} \quad \left| A_+(R, C)\right| \ \geq \ \left| A_+\left(R', C'\right) \right| \tag4.2.3$$
provided $R'$ dominates $R$ and $C'$ dominates $C$, see Chapter 16 of \cite{LW01} and \cite{Ba07}. 
Inequalities (4.2.3) are consistent with the inequalities of (4.2.1) and (4.2.2).

\head 5. Dependence between row and column sums \endhead 

The following attractive ``independence heuristic'' for estimating $\left|A_0(R, C)\right|$ and 
$\left|A_+(R, C)\right|$ was discussed by Good \cite{Go77} and by Good and 
Crook \cite{GC76}.

\subhead (5.1) The independence heuristic \endsubhead Let us consider the set of all $m \times n$ matrices 
$D=\left(d_{ij}\right)$ with 0-1 entries and the total sum $N$ of entries as a finite probability space with the 
uniform measure. Let us consider the event ${\Cal R}_0$ consisting of the matrices with the row sums 
$R=\left(r_1, \ldots, r_m \right)$ and the event ${\Cal C}_0$ consisting of the matrices with the column sums 
$C=\left(c_1, \ldots, c_n\right)$. Then 
$$\PP\left({\Cal R}_0\right) = {mn \choose N}^{-1} \prod_{i=1}^m {n \choose r_i} \quad 
\text{and} \quad \PP\left({\Cal C_0}\right)={mn \choose N}^{-1} \prod_{j=1}^n {m \choose c_j}.$$
In addition,
$$A_0(R, C) ={\Cal R}_0 \cap {\Cal C}_0.$$
If we assume that events ${\Cal R}_0$ and ${\Cal C}_0$ are independent, we obtain the following 
{\it independence estimate} 
$$I_0(R, C) ={mn \choose N}^{-1}\prod_{i=1}^m {n \choose r_i} \prod_{j=1}^n {m \choose c_j} \tag5.1.1$$
for the number $\left| A_0(R, C)\right|$ of 0-1 matrices with row sums $R$ and column sums $C$.

Similarly, let us consider the set of all $m \times n$ matrices $D=\left(d_{ij}\right)$ with non-negative integer entries and 
the total sum $N$ of entries as a finite probability space with the uniform measure. Let us consider the event 
${\Cal R}_+$ consisting of the matrices with the row sums $R=\left(r_1, \ldots, r_m \right)$ and the event 
${\Cal C}_+$ consisting of the matrices with the column sums $C=\left(c_1, \ldots, c_n \right)$. Then
$$\split \PP\left({\Cal R}_+\right) =&{N +mn-1 \choose mn-1}^{-1} \prod_{i=1}^m {r_i + n-1 \choose n-1} \quad 
\text{and} \\ \PP\left({\Cal C}_+\right)=&{N +mn-1 \choose mn-1}^{-1} \prod_{j=1}^n {c_j +m-1 \choose m-1}.
\endsplit$$
We have 
$$A_+(R, C)={\Cal R}_+ \cap {\Cal C}_+.$$
If we assume that events ${\Cal R}_+$ and ${\Cal C}_+$ are independent, we obtain the {\it independence estimate}
$$I_+(R, C) ={N +mn-1 \choose mn-1}^{-1} \prod_{i=1}^m {r_i +n-1 \choose n-1} 
\prod_{j=1}^n {c_j +m-1 \choose m-1}. \tag5.1.2$$
Interestingly, the independence estimates $I_0(R, C)$ and  $I_+(R, C)$ provide reasonable approximations 
to $\left| A_0(R, C)\right|$ and $\left| A_+(R, C)\right|$ respectively in the following two cases:
\smallskip
in the case of equal margins, when 
$$r_1 =\ldots =r_m =r \quad \text{and} \quad c_1 =\ldots = c_n =c,$$ 
see \cite{C+08} and \cite{C+07}
\smallskip
in the sparse case, when 
$$\max_{i=1, \ldots, m} r_i \ll n \quad \text{and} \quad \max_{j=1, \ldots, n} c_j \ll m,$$
see \cite{G+06} and \cite{GM08}.
\smallskip
We will see in Section 5.4 that the independence estimates provide the correct logarithmic asymptotics in the case
when all row sums are equal {\it or} all column sums are equal.
However, if both row and column sums are sufficiently far away from being uniform and sparse, the independence 
estimates, generally speaking, provide poor approximations. Moreover, in the case of 0-1 matrices the independence 
estimate $I_0(R, C)$ typically grossly overestimates $\left| A_0(R, C)\right|$ while in the case of non-negative 
integer matrices the independence estimate $I_+(R, C)$ typically grossly underestimates $\left| A_+(R, C)\right|$.
In other words, for typical margins $R$ and $C$ the events ${\Cal R}_0$ and 
${\Cal C}_0$ {\it repel} each other (the events are negatively correlated) while events ${\Cal R}_+$ and ${\Cal C}_+$ 
{\it attract} each other (the events are positively correlated). To see why this is the case, we write the estimates 
$\alpha_0(R, C)$ of Theorem 2.1 and $\alpha_+(R, C)$ of Theorem 3.1 in terms of entropy.

The following result is proven in \cite{Ba10a}. 
\proclaim{(5.2) Lemma} Let $P_0(R, C)$ be the polytope of all $m \times n$ matrices 
$X=\left(x_{ij}\right)$ with row sums $R$, column sums $C$ and such that 
$0 \leq x_{ij} \leq 1$ for all $i$ and $j$. Suppose that polytope $P_0(R, C)$ has a 
non-empty interior, that is contains a matrix $Y=\left(y_{ij}\right)$ such that 
$0 < y_{ij} < 1$ for all $i$ and $j$.  Let us define a function 
$h: P_0(R, C) \longrightarrow {\Bbb R}$ by 
$$h(X)=\sum_{i, j} x_{ij} \ln {1 \over x_{i, j}} + \left(1-x_{ij} \right) \ln {1 \over 1-x_{ij}}
\quad \text{for} \quad X \in P_0(R, C).$$
Then $h$ is a strictly concave function on of $P_0(R, C)$ and hence attains its maximum
on $P(R, C)$ at a unique matrix $Z_0=\left(z_{ij}\right)$, which we call the maximum 
entropy matrix. Moreover,
\roster 
\item We have $0 < z_{ij} < 1$ for all $i$ and $j$;
\item The infimum $\alpha_0(R, C)$ of Theorem 2.1 is attained at some particular 
point $(\xx, \yy)$;
\item We have $\alpha_0(R, C)=e^{h(Z_0)}$.
\endroster
\endproclaim
\demo{Sketch of Proof} It is straightforward to check that $h$ is strictly concave and that 
$${\partial \over \partial x_{ij}} h(X) =\ln {1-x_{ij} \over x_{ij}}.$$
In particular, the (right) derivative at $x_{ij}=0$ is $+\infty$, the (left) derivative at $x_{ij}=1$ is 
$-\infty$ and the derivative for $0 < x_{ij} < 1$ is finite.  Hence the maximum entropy 
matrix $Z_0$ must have all entries strictly between 0 and 1, since otherwise we 
can increase the value of $h$ by perturbing $Z_0$ in the direction of a matrix $Y$ 
from the interior of $P_0(R, C)$. This proves Part (1).

The Lagrange optimality conditions imply that 
$$\ln {1 -z_{ij} \over z_{ij}} =-\lambda_i - \mu_j \quad \text{for all} \quad i,j$$
and some numbers $\lambda_1, \ldots, \lambda_m$ and $\mu_1, \ldots, \mu_n$.
Hence 
$$z_{ij} = {e^{\lambda_i + \mu_j} \over 1+ e^{\lambda_i +\mu_j}} \quad \text{for all} \quad 
i, j. \tag5.2.1$$
In particular,
$$\split &\sum_{i=1}^m {e^{\lambda_i +\mu_j} \over 1+e^{\lambda_i +\mu_j}} =c_j \quad 
\text{for} \quad j=1, \ldots, n \quad \text{and} \\
&\sum_{j=1}^n {e^{\lambda_i +\mu_j} \over 1+e^{\lambda_i +\mu_j}} =r_i 
\quad \text{for} \quad i=1, \ldots, m. \endsplit \tag5.2.2$$
Equations (5.2.2) imply that point $\sss=\left(\lambda_1, \ldots, \lambda_m\right)$ and 
$\ttt=\left(\mu_1, \ldots, \mu_n \right)$ is a critical point of function $G_0(\sss, \ttt)$ defined by 
(2.2.1) and hence the infimum $\alpha_0(R, C)$ of $F_0(\xx, \yy)$ is attained at 
$x_i=e^{\lambda_i}$ for $i=1, \ldots, m$ and $y_j =e^{\mu_j}$ for $j=1, \ldots, n$. Hence 
Part (2) follows. Using (5.2.1) it is then straightforward to check that $F_0(\xx, \yy)=e^{h(Z_0)}$ 
for the minimum point $(\xx, \yy)$.
{\hfill \hfill \hfill} \qed
\enddemo

We note that 
$$h(x) =x\ln {1 \over x} +(1-x) \ln {1 \over x} \quad \text{for} \quad 0 \leq x\leq 1$$
is the entropy of the Bernoulli random variable with expectation $x$, see Section 6.

The following result is proven in \cite{Ba09}. 

\proclaim{(5.3) Lemma} Let $P_+(R, C)$ be the polytope of all non-negative $m \times n$ 
matrices $X=\left(x_{ij}\right)$ with row sums $R$ and column sums $C$. Let us define 
a function $g: P_+(R, C) \longrightarrow {\Bbb R}$ by 
$$g(X)=\sum_{i,j} \left(x_{ij}+1\right) \ln \left(1 + x_{ij}\right) -x_{ij} \ln x_{ij} 
\quad \text{for} \quad X \in P_+(R, C).$$
Then  $g$ is a strictly concave function on $P_+(R, C)$ and hence attains its maximum 
on $P_+(R, C)$ at a unique matrix $Z_+=\left(z_{ij}\right)$, which we call the maximum 
entropy matrix. Moreover,
\roster
\item We have $z_{ij} >0$ for all $i, j$ and
\item For the minimum $\alpha_+(R, C)$ of Theorem 3.1, we have 
$\alpha_+(R, C)=e^{g(Z_+)}$.
\endroster
\endproclaim 
\demo{Sketch of Proof} It is straightforward to check that $g$ is strictly concave and that 
$${\partial \over \partial x_{ij}} g(X)=\ln {1 + x_{ij} \over x_{ij}} \quad \text{for all} \quad i,j.$$
In particular, the (left) derivative is $+\infty$ for $x_{ij}=0$ and finite for every $x_{ij}>0$.
Since $P_+(R, C)$ contains an interior point (for example, matrix $Y=\left(y_{ij}\right)$ 
with $y_{ij}=r_i c_j/N$), arguing as in the proof of Lemma 5.2, we obtain Part (1).

The Lagrange optimality conditions imply that 
$$\ln {1 +z_{ij} \over z_{ij}} = \lambda_i + \mu_j \quad \text{for all} \quad i,j$$
and some numbers $\lambda_1, \ldots, \lambda_m$ and $\mu_1, \ldots, \mu_n$.
Hence 
$$z_{ij} ={e^{-\lambda_i -\mu_j} \over 1-e^{-\lambda_i -\mu_j}} \quad \text{for all} \quad i, j.
\tag5.3.1$$
In particular,
$$\split &\sum_{i=m}^n {e^{-\lambda_i -\mu_j} \over 1+ e^{-\lambda_i -\mu_j}}=c_j \quad 
\text{for} \quad j=1, \ldots, n \quad \text{and} \\
&\sum_{j=1}^n {e^{-\lambda_i -\mu_j} \over 1+e^{-\lambda_i -\mu_j}} =r_i 
\quad \text{for} \quad i=1, \ldots, m. \endsplit \tag5.3.2$$
Equations (5.3.2) imply that the point $\sss=\left(\lambda_1, \ldots, \lambda_m \right)$ 
and $\ttt=\left(\mu_1, \ldots, \mu_n\right)$ is a critical point of function $G_+(\sss, \ttt)$ defined by (3.2.1)
and hence the minimum $\alpha_+(R, C)$ of $F_+(\xx, \yy)$ is attained at 
$x_i =e^{\lambda_i}$ for $i=1, \ldots, m$ and $y_j=e^{\mu_j}$ for $j=1, \ldots, n$.
Using (5.3.1), it is then straightforward to check that $F_+(\xx, \yy)=e^{h(Z_+)}$ for the 
minimum point $(\xx, \yy)$.
{\hfill \hfill \hfill} \qed
\enddemo

We note that 
$$g(x)=(x+1) \ln (x+1) -x \ln x \quad \text{for} \quad x \geq 0$$
is the entropy of the geometric random variable with expectation $x$, see Section 6.

\subhead (5.4) Estimates of the cardinality via entropy \endsubhead
Let 
$$\HH\left(p_1, \ldots, p_k \right)=\sum_{i=1}^k p_i \ln {1 \over p_i}$$
be the entropy function defined on $k$-tuples (probability distributions) $p_1, \ldots, p_k$ such that $p_1 + \ldots + p_k=1$ 
and $p_i  \geq 0$ for $i=1, \ldots, k$. Assuming that polytope $P_0(R, C)$ of Lemma 5.2 has a non-empty 
interior, we can write
$$\aligned \ln \alpha_0(R, C)=&N \HH\left( {z_{ij} \over N }; \ i,j \right) + (mn-N) \HH \left( {1-z_{ij} \over mn -N}; \ i,j \right) \\
&\qquad -N \ln N - (mn-N) \ln (mn-N), \endaligned $$
where $Z_0=\left(z_{ij}\right)$ is the maximum entropy matrix.
On the other hand, for the independence estimate (5.1.1), we have 
$$\aligned \ln I_0(R, C)=&N \HH\left({r_i \over N};\ i\right)+(mn-N) \HH\left({n-r_i \over mn -N};\ i\right)\\
&\qquad +N \HH\left({c_j \over N};\ j \right) +(mn-N)\HH\left({m-c_j \over mn-N};\ j\right) \\
&\qquad -N \ln N -(mn-N) \ln (mn-N)+
O\bigl((m+n)\ln (mn)\bigr). \endaligned$$
Using the inequality which relates the entropy of a distribution and the entropy of its margins, see, 
for example, \cite{Kh57}, we obtain 
$$\HH \left({z_{ij} \over N};\ i, j \right) \ \leq \ \HH\left({r_i \over N};\ i\right) +\HH\left({c_j \over N}; \ j \right) \tag5.4.1$$
with the equality if and only if 
$$z_{ij} ={r_i c_j \over N} \quad \text{for all} \quad i,j$$
and 
$$\HH \left({1-z_{ij} \over mn-N};\ i, j \right) \ \leq \ \HH\left({n-r_i \over mn- N};\ i\right) +\HH\left({m-c_j \over mn-N}; \ j \right) \tag5.4.2$$
with the equality if and only if
$$1-z_{ij}={\left(n-r_i \right)\left(m-c_j\right) \over mn-N} \quad \text{for all} \quad i, j.$$
Thus we have equalities in (5.4.1) and (5.4.2) if and only if
$$\left(r_i m -N\right)\left(c_j n-N\right)=0 \quad \text{for all} \quad i, j,$$
that is, when all row sums are equal or all column sums are equal. In that case $I_0(R, C)$ 
estimates $\left|A_0(R, C)\right|$ within an $(mn)^{O(m+n)}$ factor. In all other cases, 
$I_0(R, C)$ overestimates $\left| A_0(R, C)\right|$ by as much as a $2^{\Omega(mn)}$ factor as long as the differences 
between the right hand sides and left hand sides of (5.4.1) and (5.4.2) multiplied by $N$ and $(mn-N)$ respectively
 overcome the 
$O\bigl((m+n) \ln (mn) \bigr)$ error term, see also Section 5.5 for a particular family of examples.

We handle non-negative integer matrices slightly differently. For the independence estimate (5.1.2) we obtain
$$\aligned \ln I_+(R, C)=&-(N+mn) \HH\left({r_i+n \over N+mn};\ i\right)-(N+mn) \HH\left({c_j +m \over N+mn};\ j\right)\\
&\qquad -\sum_{i=1}^m r_i \ln r_i -\sum_{j=1}^n c_j \ln c_j\\
&\qquad +N \ln N +(N+mn) \ln (N+mn)+
O\bigl((m+n)\ln N\bigr) \endaligned$$
On the other hand, by Lemma 5.3 we have 
$$\ln \alpha_+(R, C) =g\left(Z_+\right) \ \geq \ g(Y),$$
where $Z_+$ is the maximum entropy matrix and $Y=\left(y_{ij}\right)$ is the matrix defined by 
$$y_{ij}={r_i c_j \over N} \quad \text{for all} \quad i, j.$$
It is then easy to check that 
$$\split g(Y)=&-(N+mn) \HH\left( {r_i c_j +N \over N(N+mn)};\ i, j \right) \\
&\qquad -\sum_{i=1}^m r_i \ln r_i -\sum_{j=1}^n c_j \ln c_j\\
&\qquad +N \ln N +(N+mn) \ln (N+mn). \endsplit$$
By the inequality relating the entropy of a distribution and the entropy of its margins \cite{Kh57}, we have 
$$\HH\left( {r_i c_j +N \over N(N+mn)};\ i, j \right) \ \leq \ \HH\left({r_i+n \over N+mn};\ i\right) +
 \HH\left({c_j +m \over N+mn};\ j\right) \tag 5.4.3$$
with the equality if and only if 
$${r_i c_j +N \over N(N+mn)} ={(r_i +n)(c_j+m) \over (N+mn)^2} \quad \text{for all} \quad i, j,$$
that is, when we have 
$$\left(r_i m -N\right)\left(c_j n-N\right)=0 \quad \text{for all} \quad i, j,$$
so that all row sums are equal or all column sums are equal. In that case, by symmetry we have $Y=Z_+$ and 
hence $I_+(R, C)$ estimates $\left| A_+(R, C)\right|$ within an $N^{O(m+n)}$ factor. In all other cases, 
$I_+(R, C)$ underestimates $\left| A_+(R, C)\right|$ by as much as a $2^{\Omega(mn)}$ factor as long as the difference
between the right hand side and left hand side of (5.4.3) multiplied by $N+mn$ overcomes the 
$O\bigl((m+n) \ln N \bigr)$ error term, see also Section 5.5 for a particular family of examples.

\subhead (5.5) Cloning margins \endsubhead 

Let us choose a positive integer $m$-vector $R=\left(r_1, \ldots, r_m \right)$ and a positive integer $n$-vector
$C=\left(c_1, \ldots, c_n \right)$ such that 
$$r_1 +\ldots + r_m =c_1 +\ldots +c_n =N.$$
For a positive integer $k$, let us define a $km$-vector $R_k$ and a $kn$-vector $C_k$ by
$$\split R_k=&\left(\underbrace{kr_1, \ldots, kr_1}_{\text{$k$ times}}, \ldots, \underbrace{kr_m, \ldots, kr_m}_{\text{$k$ times}}\right)
\quad \text{and} \\ C_k=&\left(\underbrace{kc_1, \ldots, kc_1}_{\text{$k$ times}}, \ldots, 
\underbrace{kc_n, \ldots, kc_n}_{\text{$k$ times}}\right). \endsplit$$
We say that margins $\left(R_k, C_k\right)$ are obtained by {\it cloning} from margins $(R, C)$.
It is not hard to show that if $Z_0$ and $Z_+$ are the maximum entropy matrices associated with margins $(R, C)$ 
via Lemma 5.2 and Lemma 5.3 respectively, then the maximum entropy matrices associated with margins 
$\left(R_k, C_k \right)$ are the Kronecker products $Z_0 \otimes Id_k$ and $Z_+ \otimes Id_k$ respectively, where $Id_k$ is the $k \times k$ 
identity matrix. One has
$$\split &\lim_{k \longrightarrow +\infty} \left| A_0\left(R_k, C_k \right)\right|^{1/k^2}  =\alpha_0(R, C) \quad 
\text{and} \\ &\lim_{k \longrightarrow +\infty} \left| A_+\left(R_k, C_k\right) \right|^{1/k^2} =\alpha_+(R, C). \endsplit$$
Moreover, if not all coordinates $r_i$ of $R$ are equal and not all coordinates $c_j$ of $C$ are equal then 
the independence estimate $I_0\left(R_k, C_k\right)$, see (5.1.1), overestimates the number of $km \times kn$ matrices
with row sums $R_k$ and column sums $C_k$ and 0-1 entries within a $2^{\Omega(k^2)}$ factor while 
the independence estimate $I_+\left(R_k, C_k\right)$, see (5.1.2), underestimates the number of 
$km \times kn$ non-negative integer matrices within a $2^{\Omega(k^2)}$ factor, see \cite{Ba10a} and \cite{Ba09}
for details.

\head 6. Random matrices with prescribed row and column sums \endhead

Estimates of Theorems 2.1 and 3.1, however crude, allow us to obtain a description of a random or typical matrix 
from sets $A_0(R, C)$ and $A_+(R, C)$, considered as finite probability spaces with the uniform measures.

Recall that $x$ is a {\it Bernoulli} random variable if 
$$\PP\{x=0\}=p \quad \text{and} \quad \PP\{x=1\} =q$$
for some $p, q\geq 0$ such that $p+q=1$. Clearly, $\EE x=q$.

Recall that $P_0(R, C)$ is the polytope of $m \times n$ matrices with row sums $R$, column sums $C$ and 
entries between 0 and 1. Let function $h: P_0(R, C) \longrightarrow {\Bbb R}$ and the 
maximum entropy matrix $Z_0 \in P_0(R, C)$ be defined as in Lemma 5.2.

The following result is proven in \cite{Ba10a}, see also \cite{BH10a}.
\proclaim{(6.1) Theorem} Suppose that polytope $P_0(R, C)$ has a non-empty interior and let $Z_0 \in P_0(R, C)$ 
be the maximum entropy matrix. Let $X=\left(x_{ij}\right)$ be a random $m \times n$ 
matrix of independent Bernoulli random variables $x_{ij}$ such that $\EE X=Z_0$. 
Then 
\roster
\item The probability mass function of $X$ is constant on the set $A_0(R, C)$ of 0-1 matrices with row sums $R$ and 
column sums $C$ and 
$$\PP\{X=D\} =e^{-h\left(Z_0\right)} \quad \text{for all} \quad D \in A_0(R, C);$$
\item We have 
$$\PP\left\{X \in A_0(R, C) \right\} \ \geq \ (mn)^{-\gamma(m+n)},$$
where $\gamma >0$ is an absolute constant.
\endroster
\endproclaim

Theorem 6.1 implies that in many respects a random matrix $D \in A_0(R, C)$ behaves as a random matrix $X$ 
of independent Bernoulli random variables such that $\EE X=Z_0$, where $Z_0$ is the maximum entropy matrix.
More precisely, any event that is sufficiently rare for the random matrix $X$ (that is, an event the probability of which is 
essentially smaller than $(mn)^{-O(m+n)}$), will also be a rare event for a random matrix $D \in A_0(R, C)$. In particular,
we can conclude that a typical matrix $D \in A_0(R, C)$ is sufficiently close to $Z_0$ as long as sums of entries 
over sufficiently large subsets $S$ of indices are concerned.

For an $m \times n$ matrix $B=\left(b_{ij}\right)$ and a subset 
$$S \subset \Bigl\{ (i, j): \quad i=1, \ldots, m, \quad j =1, \ldots, n \Bigr\}$$
let 
$$\sigma_S(B)=\sum_{(i, j) \in S} b_{ij}$$
be the sum of the entries of $B$ indexed by set $S$. We obtain the following corollary, see \cite{Ba10a} for details. 
\proclaim{(6.2) Corollary} Let us fix real numbers $\kappa >0$ and $0 < \delta <1$. Then there exists a number 
$q=q(\kappa, \delta)>0$ such that the following holds.

Let $(R, C)$ be margins such that $n \geq m >q$ and the polytope $P_0(R, C)$ has a non-empty interior and let 
$Z_0 \in P_0(R, C)$ be the maximum entropy matrix. Let $S \subset \bigl\{(i, j): \quad i=1, \ldots, m; \ 
j=1, \ldots, n \bigr\}$ be a set such that $\sigma_S\left(Z_0\right) \geq \delta mn$ and let 
$$\epsilon =\delta {\ln \over \sqrt{m}}.$$
If $\epsilon \leq 1$ then
$$\PP \Bigl\{ D \in A_0(R, C): \ (1-\epsilon) \sigma_S\left(Z_0\right) \ \leq \ \sigma_S(D) \ \leq \ 
(1+\epsilon) \sigma_S\left(Z_0\right) \Bigr\} \ \geq \ 1 -n^{-\kappa n}.$$
\endproclaim

Recall that $x$ is a {\it geometric} random variable if 
$$\PP\{x=k\}=pq^k \quad \text{for} \quad k=0, 1, 2, \ldots $$
for some $p, q\geq 0$ such that $p+q=1$. We have $\EE x=q/p$.

Recall that $P_+(R, C)$ is the polytope of $m \times n$ non-negative matrices with row sums $R$ and column sums $C$.
Let function $g: P_+(R, C) \longrightarrow {\Bbb R}$ and the 
maximum entropy matrix $Z_+ \in P_0(R, C)$ be defined as in Lemma 5.3.

The following result is proven in \cite{Ba10b}, see also \cite{BH10a}.
\proclaim{(6.3) Theorem} Let $Z_+ \in P_0(R, C)$ 
be the maximum entropy matrix. Let $X=\left(x_{ij}\right)$ be a random $m \times n$ 
matrix of independent geometric random variables $x_{ij}$ such that $\EE X=Z_+$. 
Then 
\roster
\item The probability mass function of $X$ is constant on the set $A_+(R, C)$ of non-negative
integer matrices with row sums $R$ and 
column sums $C$ and 
$$\PP\{X=D\} =e^{-g\left(Z_+\right)} \quad \text{for all} \quad D \in A_+(R, C);$$
\item We have 
$$\PP\left\{X \in A_+(R, C) \right\} \ \geq \ N^{-\gamma(m+n)},$$
where $\gamma >0$ is an absolute constant and $N=r_1 +\ldots + r_m =c_1 + \ldots + c_n$ 
for $R=\left(r_1, \ldots, r_m \right)$ and $C=\left(c_1, \ldots, c_n \right)$.
\endroster
\endproclaim

Theorem 6.3 implies that in many respects a random matrix $D \in A_+(R, C)$ behaves as a matrix $X$ 
of independent geometric random variables such that $\EE X=Z_+$, where $Z_+$ is the maximum entropy matrix.
More precisely, any event that is sufficiently rare for the random matrix $X$ (that is, an event the probability of which is 
essentially smaller than $N^{-O(m+n)}$), will also be a rare event for a random matrix $D \in A_+(R, C)$. In particular,
we can conclude that a typical matrix $D \in A_+(R, C)$ is sufficiently close to $Z_+$ as long as sums of entries 
over sufficiently large subsets $S$ of indices are concerned.

Recall that $\sigma_S(B)$ denotes the sum of the entries of a matrix $B$ indexed by a set $S$.
We obtain the following corollary, see \cite{Ba10b} for details. 

\proclaim{(6.4) Corollary} Let us fix real numbers $\kappa >0$ and $0 < \delta <1$. Then there exists a positive integer 
$q=q(\kappa, \delta)$ such that the following holds.

Let $R=\left(r_1, \ldots, r_m \right)$ and $C=\left(c_1, \ldots, c_m\right)$ be positive integer vectors such that 
$r_1+ \ldots + r_m =c_1 + \ldots +c_n =N$,
$$\split &{\delta N \over m} \ \leq \ r_i \ \leq \ {N \over \delta m} \quad \text{for} \quad i=1, \ldots m, \\
&{\delta N \over n} \ \leq c_j \ \leq \ {N \over \delta n} \quad \text{for} \quad j=1, \ldots, n \endsplit$$
and 
$${N \over mn} \ \geq \ \delta.$$
Suppose that $n \geq m >q$ and let 
$S \subset \bigl\{(i, j): \quad i=1, \ldots, m, \ j=1, \ldots, n \bigr\}$ be a set such that 
$|S| \geq \delta mn$.
Let $Z_+ \in P_+(R, C)$ be the maximum entropy matrix and let 
$$\epsilon =\delta {\ln n \over m^{1/3}}.$$
If $\epsilon \leq 1$ then 
$$\PP\Bigl\{D \in A_+(R, C): \ (1-\epsilon) \sigma_S\left(Z_+\right) \ \leq \ \sigma_S(D) \ \leq \ 
(1+\epsilon) \sigma_S\left(Z_+\right) \Bigr\} \ \geq \ 1- n^{-\kappa n}.$$
\endproclaim

As is discussed in \cite{BH10a}, the ultimate reason why Theorems 6.1 and 6.3 hold true is that
\smallskip
 the matrix $X$ of 
independent Bernoulli random variables such that $\EE X=Z_0$ is the random matrix with the maximum 
possible entropy among all random $m \times n$ matrices with 0-1 entries and the expectation in the affine 
subspace of the matrices with row sums $R$ and column sums $C$ 
\smallskip
and
\smallskip
 the matrix $X$ of independent geometric random variables such that $\EE X=Z_+$ is the random matrix with 
the maximum possible entropy among all random $m \times n$ matrices with non-negative integer entries and the 
expectation in the affine subspace of the matrices with row sums $R$ and column sums $C$.
\smallskip
Thus Theorems 6.1 and 6.3 can be considered as an illustration of the Good's thesis \cite{Go63} that the 
``null hypothesis'' for an unknown probability distribution from a given class should be the hypothesis that the 
unknown distribution is, in fact, the distribution of the maximum entropy in the given class.

\subhead (6.5) Sketch of proof of Theorem 6.1 \endsubhead Let $Z_0=\left(z_{ij}\right)$ be the maximum entropy matrix 
as in Lemma 5.2.
Let us choose $D \in A_0(R, C)$, $D=\left(d_{ij}\right)$. Using (5.2.1), we get
$$\split \PP\bigl\{X =D \bigr\} =&\prod_{i, j} z_{ij}^{d_{ij}} \left(1 -z_{ij}\right)^{1-d_{ij}} =
\prod_{ij} {e^{\left(\lambda_i +\mu_j \right)d_{ij}} \over 1+e^{\lambda_i +\mu_j}}
\\=&\exp\left\{ \sum_{i=1}^m \lambda_i r_i + \sum_{j=1}^n \mu_j c_j \right\} 
\prod_{ij} {1 \over 1+e^{\lambda_i +\mu_j}}\\
=&e^{-h\left(Z_0\right)}, \endsplit$$
which proves Part (1).

To prove Part (2), we use Part (1), Theorem 2.1 and Lemma 5.2. We have
$$\split \PP\bigl\{X \in A_0(R, C) \bigr\} =&\left| A_0(R, C)\right| e^{-h\left(Z_0\right)} \ \geq \ (mn)^{-\gamma(m+n)} 
\alpha_0(R, C) e^{-h\left(Z_0\right)} \\=&(mn)^{-\gamma(m+n)} \endsplit$$
for some absolute constant $\gamma >0$.
{\hfill \hfill \hfill} \qed

\subhead (6.6) Sketch of proof of Theorem 6.3 \endsubhead Let $Z_+=\left(z_{ij}\right)$ be the maximum entropy matrix 
as in Lemma 5.3.
Let us choose $D \in A_+(R, C)$, $D=\left(d_{ij}\right)$. Using (5.3.1), we get
$$\split \PP\bigl\{X =D \bigr\} =&\prod_{i, j} \left({1 \over 1+ z_{ij}}\right) \left({z_{ij} \over 1+z_{ij}}\right)^{d_{ij}} =
\prod_{ij} \left(1-e^{-\lambda_i -\mu_j}\right) e^{-\left(\lambda_i +\mu_j \right) d_{ij}}
\\=&\exp\left\{ -\sum_{i=1}^m \lambda_i r_i - \sum_{j=1}^n \mu_j c_j \right\} 
\prod_{ij} \left(1-e^{-\lambda_i -\mu_j}\right)\\
=&e^{-g\left(Z_+\right)}, \endsplit$$
which proves Part (1).

To prove Part (2), we use Part (1), Theorem 3.1 and Lemma 5.3. We have
$$\split \PP\bigl\{X \in A_+(R, C) \bigr\} =&\left| A_+(R, C)\right| e^{-g\left(Z_+\right)} \ \geq \ N^{-\gamma(m+n)} 
\alpha_+(R, C) e^{-g\left(Z_+\right)} \\=&N^{-\gamma(m+n)} \endsplit$$
for some absolute constant $\gamma >0$.
{\hfill \hfill \hfill} \qed

\subhead (6.7) Open questions \endsubhead 
Theorems 6.1 and 6.3 show that a random matrix $D \in A_0(R, C)$, respectively $D \in A_+(R, C)$, in many respects 
behaves like a matrix of independent Bernoulli, respectively geometric, random variables whose expectation is the 
maximum entropy matrix $Z_0$, respectively $Z_+$. One can ask whether individual entries $d_{ij}$ of $D$ behave 
asymptotically as Bernoulli, respectively geometric, random variables with expectations $z_{ij}$ as the size of the 
matrices grows. In the simplest situation we ask the following
\subsubhead (6.7.1) Question \endsubsubhead Let $(R, C)$ be margins and let $\left(R_k, C_k\right)$ be margins 
obtained from $(R, C)$ by cloning as in Section 5.5. Is it true that as $k$ grows, the entry $d_{11}$ of a random matrix $D \in A_0\left(R_k, C_k\right)$,
respectively $D \in A_+\left(R_k, C_k\right)$, converges in distribution to the Bernoulli, respectively geometric,
random variable with expectation $z_{11}$, where $Z_0=\left(z_{ij}\right)$, respectively $Z_+=\left(z_{ij}\right)$, is 
the maximum entropy matrix of margins $(R, C)$? 
\smallskip
Some entries of the maximum entropy matrix $Z_+$ may turn out to be surprisingly large, even for reasonably 
looking margins. 
In \cite{Ba10b}, the following example is considered. Suppose that $m=n$ and let 
$R_n=C_n=(3n, n, \ldots, n)$. It turns out that the entry $z_{11}$ of the maximum entropy matrix $Z_+$ is linear 
in $n$, namely $z_{11} > 0.58n$, while all other entries remain bounded by a constant. One can ask whether the 
$d_{11}$ entry of a random matrix $D \in A_+\left(R_n, C_n\right)$ is indeed large, as the value of $z_{11}$ suggests.

\subsubhead (6.7.2) Question \endsubsubhead Let $\left(R_n, C_n\right)$ be margins as above. Is it true that as $n$ grows, 
one has $\EE d_{11} =\Omega(n)$ for a random matrix $D \in A_+\left(R_n, C_n\right)$?
\smallskip
Curiously, the entry $z_{11}$ becomes bounded by a constant if $3n$ is replaced by $2n$.

\head 7. Asymptotic formulas for the number of matrices with prescribed row and column sums
\endhead

In this section, we discuss asymptotically exact estimates for $\left| A_0(R, C)\right|$ and $\left| A_+(R, C)\right|$.

\subhead (7.1) An asymptotic formula for $\left|A_0(R, C)\right|$ \endsubhead
Theorem 6.1 suggests the following way to estimate the number $\left| A_0(R, C)\right|$ of 
0-1 matrices with row sums $R$ and column sums $C$. Let us consider the matrix 
of independent Bernoulli random variables as in Theorem 6.1 and let 
$Y$ be the random $(m+n)$-vector obtained by computing the row and column sums of $X$.
Then, by Theorem 6.1, we have 
$$\left|A_0(R, C)\right| = e^{h\left(Z_0\right)} \PP\bigl\{X \in A_0(R, C)\bigr\} 
=e^{h\left(Z_0\right)} \PP\bigl\{ Y=(R, C) \bigr\}. \tag7.1.1$$
Now, random $(m+n)$-vector $Y$ is obtained as a sum of $mn$ independent random vectors
and $\EE Y =(R, C)$, so it is not unreasonable to assume that 
$\PP\bigl\{ Y=(R, C)\bigr\}$ can be estimated via some version of the Local Central Limit Theorem.
In \cite{BH10b} we show that this is indeed the case provided one employs the Edgeworth 
correction factor in the Central Limit Theorem.

We introduce the necessary objects to state the asymptotic formula for the number of 
0-1 matrices with row sums $R$ and column sums $C$.

Let $Z_0=\left(z_{ij}\right)$ be the maximum entropy matrix as in Lemma 5.2. 
We assume that $0 < z_{ij} < 1$ for all $i$ and $j$.
Let us consider the quadratic form $q_0: {\Bbb R}^{m+n} \longrightarrow {\Bbb R}$ defined by 
$$\split q_0(s, t)=&{1 \over 2} \sum \Sb 1 \leq i \leq m \\ 1 \leq j \leq n \endSb \left(z_{ij}-z_{ij}^2\right)
\left(s_i +t_j \right)^2 \\ &\qquad \text{for} \quad s=\left(s_1, \ldots, s_m \right) \quad \text{and} 
\quad t=\left(t_1, \ldots, t_n \right). \endsplit$$
Quadratic form $q_0$ is positive semidefinite with the kernel spanned by vector 
$$u=\left(\underbrace{1, \ldots, 1}_{\text{$m$ times}}; \underbrace{-1, \ldots, -1}_{\text{$n$ times}}\right).$$
Let $H=u^{\bot}$ be the hyperplane in ${\Bbb R}^{m+n}$ defined by the equation 
$$s_1 + \ldots  + s_m = t_1 + \ldots + t_n. \tag7.1.2$$
Then the restriction $q_0|H$ of $q_0$ onto $H$ is a positive definite quadratic form 
and we define its determinant $\det q_0|H$ as the product of the non-zero eigenvalues of 
$q_0$. We consider the Gaussian probability measure on $H$ with the density proportional to 
$e^{-q_0}$ and define random variables $\phi_0, \psi_0: H\longrightarrow {\Bbb R}$ by  
$$\split \phi_0(s, t)=&{1 \over 6} \sum \Sb 1 \leq i \leq m \\ 1 \leq j \leq n \endSb 
z_{ij}\left(1-z_{ij}\right)\left(2 z_{ij} -1 \right) \left(s_i + t_j \right)^3 \quad \text{and} \\
\psi_0(s, t)=&{1 \over 24} \sum \Sb 1 \leq i \leq m \\ 1 \leq j \leq n \endSb 
z_{ij}\left(1-z_{ij}\right) \left(6z_{ij}^2 - 6z_{ij}+1 \right) \left(s_i + t_j \right)^4 \\
&\qquad \text{for} \quad (s, t)=\left(s_1, \ldots, s_m; t_1, \ldots, t_n \right). \endsplit$$
We let 
$$\mu_0 =\EE \phi_0^2 \quad \text{and} \quad \nu_0 =\EE \psi_0.$$

\proclaim{(7.2) Theorem} Let us fix $0 < \delta < 1/2$, let $R=\left(r_1, \ldots, r_m\right)$ 
and $C=\left(c_1, \ldots, c_n \right)$ be margins such that $m \geq \delta n$ and $n \geq \delta m$.
Let $Z_0=\left(z_{ij}\right)$ be the maximum entropy matrix as in Lemma 5.2 and suppose that 
$\delta \leq z_{ij} \leq 1-\delta$ for all $i$ and $j$. 

Let the quadratic form $q_0$ and values $\mu_0$ and $\nu_0$ be as defined in Section 7.1.
Then the number 
$${e^{h\left(Z_0\right)} \sqrt{m+n} \over (4 \pi)^{(m+n-1)/2} \sqrt{\det q_0|H}} 
\exp\left\{ -{\mu_0 \over 2} +\nu_0 \right\} \tag7.2.1$$
approximates the number $\left| A_0(R, C)\right|$ of 
as $m, n \longrightarrow +\infty$ within a relative error which approaches 0 as $m, n \longrightarrow +\infty$. More precisely, for any $0 < \epsilon \leq 1/2$, the value of 
(7.2.1) approximates $\left|A_0(R, C)\right|$ within relative error $\epsilon$ provided 
$$m, n\  \geq \ \left({1 \over \epsilon}\right)^{\gamma(\delta)}$$
for some $\gamma(\delta) >0$.
\endproclaim

Some remarks are in order.

All the ingredients of formula (7.2.1) are efficiently computable, in time polynomial in $m+n$, see \cite{BH10b} for details.
If all row sums are equal then we have $z_{ij}=c_j/m$ by symmetry and if all column sums are equal, we have 
$z_{ij}=r_i/n$. In particular, if all row sums are equal and if all column sums are equal, we obtain the asymptotic 
formula of \cite{C+08}.

Let us consider formula (7.1.1). If, in the spirit of the Local Central Limit Theorem, we approximated
$\PP\bigl\{Y =(R, C)\bigr\}$ by $\PP\bigl\{Y^{\ast} \in (R, C) + \Pi \bigr\}$, where $Y^{\ast}$ is the 
$(m+n-1)$-dimensional random Gaussian vector whose expectation and covariance matrix match those of $Y$ 
and where $\Pi$ is the set of points on the hyperplane $H$ that are closer to $(R, C)$ than to any other integer vector 
in $H$, we would have obtained the first part 
$${e^{h\left(Z_0\right)} \sqrt{m+n} \over (4 \pi)^{(m+n-1)/2} \sqrt{\det q_0|H}} $$
of formula (7.2.1). Under the conditions of Theorem 7.2 we have 
$$c_1(\delta) \ \leq \ \exp\left\{ -{\mu_0 \over 2} + \nu_0 \right\} \ \leq \ c_2(\delta)$$ 
for some constants $c_1(\delta), c_2(\delta) >0$ and this factor represents the Edgeworth correction 
to the Central Limit Theorem.  We note that the constraints $\delta \leq z_{ij} \leq 1-\delta$ are, generally speaking, 
unavoidable. If the entries $z_{ij}$ of the maximum entropy matrix are uniformly small, then the distribution of the random vector $Y$ of row and column sums of the random Bernoulli matrix $X$ is no longer approximately Gaussian
but approximately Poisson and formula (7.2.1) does not give correct asymptotics. The sparse case of small row and column sums is
investigated in \cite{G+06}.

More generally, to have some analytic formula approximating $\left| A_0(R, C)\right|$ we need certain regularity 
conditions on $(R, C)$, since the number $\left| A_0(R, C)\right|$ becomes volatile when the margins $(R, C)$ approach 
the boundary of the Gale-Ryser conditions, cf. \cite{JSM92}. By requiring that the entries of maximum entropy matrix $Z_0$
are separated from both 0 and 1, we ensure that the margins $(R, C)$ remain sufficiently inside the 
polyhedron defined by the Gale-Ryser inequality and the number of 0-1 matrices with row sums $R$ and 
column sums $C$ changes sufficiently smoothly when $R$ and $C$ change.

\subhead (7.3) An asymptotic formula for $\left| A_+(R, C)\right|$ \endsubhead
As in Theorem 6.3, let $X$ be the matrix of independent geometric random variables 
such that $\EE X=Z_+$, where $Z_+$ is the maximum entropy matrix. Let $Y$ be the random $(m+n)$-vector obtained by computing the row and column sums of $X$. Then, by Theorem 6.3, we have 
$$\left|A_+(R, C)\right|=e^{g\left(Z_+\right)} \PP\bigl\{X \in A_+(R, C)\bigr\} =
e^{g(Z_+)} \PP\bigl\{Y =(R, C)\bigr\}. \tag7.3.1$$
In \cite{BH09} we show how to estimate the probability that $Y=(R, C)$ using the 
Local Central Limit Theorem with the Edgeworth correction. 

Let $Z_+=\left(z_{ij}\right)$ be the maximum entropy matrix as in Lemma 5.3. Let us consider 
the quadratic form $q_+: {\Bbb R}^{m+n} \longrightarrow {\Bbb R}$ defined by 
$$\split q_+(s, t)=&{1 \over 2} \sum \Sb 1 \leq i \leq m \\ 1 \leq j \leq n \endSb \left(z_{ij} +z_{ij}^2\right)
\left(s_i +t_j \right)^2  \\
&\qquad \text{for} \quad s=\left(s_1, \ldots, s_m \right) \quad \text{and} \quad 
t=\left(t_1, \ldots, t_n \right). \endsplit$$
Let $H \subset {\Bbb R}^{m+n}$ be the hyperplane defined by (7.1.2). The restriction 
$q_+|H$ of $q_+$ onto $H$ is a positive definite quadratic form and we define its determinant 
$\det q_+|H$ as the product of the non-zero eigenvalues of $q_+$. We consider the Gaussian 
probability measure on $H$ with the density proportional to $e^{-q_+}$ and define random 
variables $\phi_+, \psi_+: H \longrightarrow {\Bbb R}$ by 
$$\split \phi_+(s, t)=&{1 \over 6} \sum \Sb 1 \leq i \leq m \\ 1 \leq j \leq n \endSb 
z_{ij}\left(1+z_{ij}\right)\left(2 z_{ij} +1 \right) \left(s_i + t_j \right)^3 \quad \text{and} \\
\psi_+(s, t)=&{1 \over 24} \sum \Sb 1 \leq i \leq m \\ 1 \leq j \leq n \endSb 
z_{ij}\left(1+z_{ij}\right) \left(6z_{ij}^2 + 6z_{ij}+1 \right) \left(s_i + t_j \right)^4 \\
&\qquad \text{for} \quad (s, t)=\left(s_1, \ldots, s_m; t_1, \ldots, t_n \right). \endsplit$$
We let 
$$\mu_+ =\EE \phi_+^2 \quad \text{and} \quad \nu_+ =\EE \psi_+.$$
\proclaim{(7.4) Theorem} Let us fix $0 < \delta < 1$, let $R=\left(r_1, \ldots, r_m\right)$ 
and $C=\left(c_1, \ldots, c_n \right)$ be margins such that $m \geq \delta n$ and $n \geq \delta m$.
Let $Z_+=\left(z_{ij}\right)$ be the maximum entropy matrix as in Lemma 5.3. Suppose that 
$$\delta \tau \ \leq \ z_{ij} \ \leq \ \tau \quad \text{for all} \quad i, j$$ 
for some $\tau \geq \delta$.

Let the quadratic form $q_+$ and values $\mu_+$ and $\nu_+$ be as defined in Section 7.3.
Then the number 
$${e^{g\left(Z_+\right)} \sqrt{m+n} \over (4 \pi)^{(m+n-1)/2} \sqrt{\det q_+|H}} 
\exp\left\{ -{\mu_+ \over 2} +\nu_+ \right\} \tag7.4.1$$
approximates the number $\left| A_+(R, C)\right|$ of 
as $m, n \longrightarrow +\infty$ within a relative error which approaches 0 as $m, n \longrightarrow +\infty$. More precisely, for any $0 < \epsilon \leq 1/2$, the value of 
(7.4.1) approximates $\left|A_+(R, C)\right|$ within relative error $\epsilon$ provided 
$$m, n\  \geq \ \left({1 \over \epsilon}\right)^{\gamma(\delta)}$$
for some $\gamma(\delta)>0$.
\endproclaim

All the ingredients of formula (7.4.1) are efficiently computable, in time polynomial in $m+n$, see \cite{BH09} for details.
If all row sums are equal then we have $z_{ij}=c_j/m$ by symmetry and if all column sums are equal, we have 
$z_{ij}=r_i/n$. In particular, if all row sums are equal and if all column sums are equal, we obtain the asymptotic 
formula of \cite{C+07}. The term 
$${e^{g\left(Z_+\right)} \sqrt{m+n} \over (4 \pi)^{(m+n-1)/2} \sqrt{\det q_+|H}} $$
corresponds to the Gaussian approximation for the distribution of the random vector $Y$ in (7.3.1), while 
$$\exp\left\{ -{\mu_+ \over 2} +\nu_+ \right\}$$
is the Edgeworth correction factor.

While the requirement that the entries of the maximum entropy matrix $Z_+$ are separated from 
0 is unavoidable (if $z_{ij}$ are small, the coordinates of $Y$ are asymptotically Poisson, not 
Gaussian, see \cite{GM08} for the analysis of the sparse case), it is not clear whether the 
requirement that all $z_{ij}$ are within a constant factor of each other is indeed needed. 
It could be that around certain margins $(R, C)$ the number $\left|A_+(R, C)\right|$ experiences
sudden jumps,
as the margins change, which precludes the existence of an analytic expression similar to 
(7.4.1) for $\left|A_+(R, C)\right|$. A candidate for such an abnormal behavior is supplied by the margins 
discussed in Section 6.7. Namely, if $m=n$ and $R=C=\left(\lambda n, n, \ldots, n \right)$ then for 
$\lambda=2$ all the entries of the maximum entropy matrix $Z_+$ are $O(1)$, while for $\lambda=3$ the 
first entry $z_{11}$ grows linearly in $n$. Hence for some particular $\lambda$ between $2$ and $3$ a 
certain ``phase transition''' occurs: the entry $z_{11}$ jumps from $O(1)$ to $\Omega(n)$. It would be interesting to find 
out if there is indeed a sharp change in $\left| A_+(R, C)\right|$ when $\lambda$ changes from 2 to 3.

\head 8. Concluding remarks \endhead

Method of Sections 6 and 7 have been applied to some related problems, such as counting higher-order
``tensors'' with 0-1 or non-negative integer entries and prescribed sums along coordinate hyperplanes \cite{BH10a}
and counting graphs with prescribed degrees of vertices \cite{BH10b}, which corresponds to counting symmetric 
0-1 matrices with zero trace and prescribed row (column) sums.

 In general, the problem can be described as follows: 
we have a polytope $P \subset {\Bbb R}^d$ defined as the intersection of the non-negative orthant ${\Bbb R}^d_+$ 
with an affine subspace ${\Cal A}$ in ${\Bbb R}^d$ and we construct a $d$-vector $X$ of independent Bernoulli (in the 0-1 case)
or geometric (in the non-negative integer case) random variables, so that the expectation of $X$ lies in ${\Cal A}$
and the distribution of $X$ is uniform, when restricted onto the set of 0-1 or integer points in $P$. 
Random vector $X$ is determined by its expectation $\EE X=z$ and $z$ is found by solving a convex optimization 
problem on $P$.
Since vector $X$ conditioned on the set of 0-1 or non-negative integer vectors in $P$ is uniform, the number of 0-1 or non-negative integer points in $P$ is expressed in terms of the 
probability that $X$ lies in ${\Cal A}$.
Assuming that 
the affine subspace ${\Cal A}$ is defined by a system $Ax=b$ of linear equations, where 
$A$ is $k \times d$ matrix of rank $k < d$, we define a $k$-vector $Y=AX$ of random variables 
and estimate the probability that $Y=b$ by using a Local Central Limit Theorem type argument. Here we essentially 
use that $\EE Y =b$, since the expectation of $X$ lies in ${\Cal A}$.

Not surprisingly, 
the argument works the easiest when the codimension $k$ of the affine subspace (and hence the dimension of vector $Y$) is small.
In particular, counting higher-order ``tensors'' is easier than counting matrices, the need in the Edgeworth 
correction factor, for example, disappears as the vector $Y$ turns out to be closer in distribution to a Gaussian vector, 
see \cite{BH10a}. Once a Gaussian or almost Gaussian estimate for the probability $\PP\bigl\{Y=b\bigr\}$ is established,
one can claim a certain concentration of a random 0-1 or integer point in $P$ around $z=\EE X$.

\enddocument
\end